\documentclass[leqno]{article}  
\usepackage{amsmath,amstext,amsthm,amssymb}   
\numberwithin{equation}{section}

\allowdisplaybreaks
\swapnumbers
\theoremstyle{plain}

\def\bthm#1.#2 #3\ethm{
\begin{\ifx#1ttheorem\fi%
\ifx#1llemma\fi% 
\ifx#1ccorollary\fi% 
\ifx#1pproposition\fi%
\ifx#1ddefinition\fi}
\label{#1.#2}  
#3 \end{\ifx#1ttheorem\fi%
\ifx#1llemma\fi%
\ifx#1ccorollary\fi% 
\ifx#1pproposition\fi%
\ifx#1ddefinition\fi}}

\def\t#1/{Theorem~\ref{t#1}}
\def\c#1/{Corollary~\ref{c#1}}
\def\l#1/{Lemma~\ref{l#1}}
\def\s#1/{Section~\ref{s#1}}
\def\e#1/{(\ref{e#1})}
\def\d#1/{Definition~\ref{d#1}}
\def\Label #1 {\label{#1}}

\def\norm#1.#2.{\lVert#1\rVert_{#2}}
\def\Norm#1.#2.{\bigl\lVert#1\bigr\rVert_{#2}}
\def\NOrm#1.#2.{\Bigl\lVert#1\Bigr\rVert_{#2}}
\def\NORm#1.#2.{\biggl\lVert#1\biggr\rVert_{#2}}
\def\NORM#1.#2.{\Biggl\lVert#1\Biggr\rVert_{#2}}

\def\ip#1,#2.{\langle #1,#2\rangle}
\def\Ip#1,#2.{\bigl\langle#1,#2\bigr\rangle}
\def\IP#1,#2.{\Bigl\langle#1,#2\Bigr\rangle}

\def\abs#1{\lvert#1\rvert}

\newcommand{\za}{\alpha}

\newcommand{\zve}{\varepsilon}
\newcommand{\zvf}{\varphi}

\newcommand{\zF}{\Phi}

\newcommand{\zI}{\infty}
\newcommand{\zl}{\lambda}

\newcommand{\zW}{\Omega}

\newcommand{\zs}{\sigma}

\newcommand{\zp}{\pi}

\def\z#1#2{\ifcase#1 {{\cal #2}}
\or {{\tilde{#2}}}
\or { {\hat{#2}}}
\or{{\check{#2}}}
\or {{\acute{#2}}}
\or {{\grave{#2}}}
\or {{\bar{#2}}}
\or {\dot{#2}}
\or {\overline{#2}}
\or {{\cal #2}}\fi}
 
\def\ZR{{\mathbb R}}
\def\ZZ{{\mathbb Z}}

\def\I{{\bf I}}
\def\S{{\bf S}}

\def\ind#1{{\mathbb I}_{#1}}

\def\md#1#2\endmd{\ifx0#1
\begin{equation*} #2 \end{equation*}\fi  %  single line display, no number
\ifx1#1\begin{equation}#2\end{equation}\fi   % single line display, number
\ifx2#1\begin{align*}#2\end{align*}\fi   % aligned display, no number
\ifx3#1\begin{align}#2\end{align}\fi    % aligned display, number
\ifx4#1\begin{gather*}#2\end{gather*}\fi  % multiline, not align, no number
\ifx5#1\begin{gather}#2\end{gather}\fi   % multinline, not align
\ifx6#1\begin{multline*}#2\end{multline*}\fi  %  display too long for one line
\ifx7#1\begin{multline}#2\end{mutline}\fi  % as above, with numbers
}

\def\S{{\bf S}}
\def\svf#1 {\varphi_{#1}}
\def\sw#1  {W_{#1}}
\def\sqw#1 {\sqrt{\abs{W_{#1}}}}
\def\szw#1 {\Omega_{#1}}

\begin{document}

\title{ A Weak--type Orthogonality Principle }
\author{Jose Barrionuevo and Michael T. Lacey\footnote{This author has been supported by an NSF grant    
 DMS--9706884} 
\\\small {University of South Alabama and Georgia Institute of Technology}}  
\date{}

\maketitle

%\abstract{We prove a weak type estimate for operators of the form 
%$  f \to \sum_{s\in\S}\langle f,\svf s \rangle \svf s $ for certain collections of Schwartz 
%functions $\{ \svf s \}_{s\in\S}$. This extends some of the orthogonality issues involved in the 
%study of the bilinear Hilbert transform by Lacey and Thiele. }

 \section{Introduction and Principal Inequalities}
 
We are interested in the relationships between three different concepts,  first and foremost 
is that
of the phase space, by which we generally  mean the Euclidean space formed from the cross 
product of the spatial
variable with the dual frequency variable.  Next, we want to associate subsets of that space 
with
functions,  the subset describing the location of the function in natural ways.  And finally, 
we
want to understand the extent to which orthogonality of the functions can be quantified by
geometric conditions on the corresponding sets in the phase plane.  

 These concerns are not currently very much in the forefront
of harmonic analysis, but rather the means towards an end.  We treat them as a subject in 
their own
right because the inequalities that we obtain are of  an optimal nature and they 
 refine  the basic orthogonality
issues in the proof of the bilinear Hilbert transform inequalities \cite{lt.pnas},
 and complement investigations into
``best basis" signal or image processing \cite{w}, including the directional issues that 
arise in
the context of brushlets \cite{brush}.

We state the principal results, and then turn to complementary issues and discussion. 
For a set $W\subset\ZR^d$ of finite volume we define $\zl W:=\{c(W)+\zl(x-c(W))\mid x\in W\}$ 
where
$c(W)$ is the center of mass of $W$.  
We say that $W$ is {\sl symmetric} if $-(W-c(W))=W-c(W)$.
 We call the product $s=W_s\times\zW_s$  of a symmetric convex set $W_s$ and a second set 
a {\sl tile}.  Here, we need not assume that the second set lies in $\ZR^d$, it could lie in 
some other set altogether.

We use tiles to study the connection between geometry of the phase plane and orthogonality 
which is
the intention of the  following definition.  

\bthm d.adapted 
 Let $\S$ be a set of tiles.  The functions ${\bf \zF}=\{\zvf_s\mid s\in\S\}$ are
said to be {\sl adapted  to} $\S$  if the 
functions $\zvf_s$ are Schwartz functions, and for all $s\in\S$,
\md5
\Label e.=1   \norm \zvf_s.2.=1,
\\
\label{e.perp} \ip \zvf_s,\zvf_{s'}.=0\quad\text{if $s'\in\S$, $\zW_s\cap 
\zW_{s'}=\emptyset$,}
\\
\label{e.space}
   \abs{\svf s (x)}\le{}\frac{C_0}{\sqw s }(1+\zs(x,\sw s ))^{-2d-5},\quad x\in\ZR^d,
\\
\noalign{\noindent where  we define $\zs(x,W)$ below.  }
\nonumber   \zs(x,W):=\inf\{a>0\mid x\in aW\}.
\endmd
\ethm

We call a collection of sets $\z0G$ a {\sl grid} if for all $G,G'\in\z0G$, we have $G\cap
G'=\emptyset, G$ or $G'$.

\bthm t.count1 
Let $\S$ be a set of tiles such that 
\md5
\Label e.grid \{\szw s \mid s\in\S\} \text{ ${}$ is a grid,}
\\
\Label e.disjoint  \{s\in\S\}\text{${}$ are pairwise disjoint,}
\\
\Label e.geo  s,s'\in\S,\ \szw s' \supset\szw s \text{\quad 
implies
\quad} \sw s' -c(\sw s' )\subset\sw s -c(\sw s ).
\endmd
Let $\{\svf s \mid s\in\S\}$ be adapted to $\S$, and for $\zl>0$ and
 $f\in L^2(\ZR^d)$ let 
\md1 \Label e.> 
\S_\zl:=\Bigl\{s\in\S\mid \frac{\abs{\ip f,\svf s .}}{\sqw s }\ge\zl
\Bigr\}.
\endmd
Then  we have the inequality
\md0
\sum_{s\in\S_\zl}\abs{\sw s }\le{}(2+KC_0^2)\zl^{-2}\lVert 
f\rVert_2^2.
\endmd
Here and throughout $K$ denotes a constant depending only on dimension $d$.
This inequality can be rephrased as 
\md0
\NOrm \Bigl\{\frac{\abs{\ip f,\svf s .}}{\sqw s }\ind {W_s} \mid s\in\S\Bigr\}
.L^{2,\zI}(\ZR^d\times \S).
\le{}\sqrt{2+KC_0^2}\norm f.2. .
\endmd
In this inequality, $\ind A$ denotes the indicator function for the set $A$, 
$L^{2,\zI}$ denotes the weak $L^2$ space, and we assign
$\ZR^d\times \S$ the product measure of Lebesque measure times counting measure
on $\S$.
\ethm

\def\C{{\bf C}}  \def\I{{\bf I}}

Note that in the last inequality if the weak $L^2$ norm could be replaced by the
 $L^2$ norm, we would have  a Bessel inequality. 
 However, 
 the weak $L^2$ space   cannot  be replaced by 
 any smaller Lorentz space.  We demonstrate this with an example at the conclusion of the paper.
 However, in the context of computation, one cannot distinguish between $L^2$ and 
weak--$L^2$.  

It is worth noting that the concluding Lemmas of \cite{lt.pnas} and \cite{lt} study 
exactly the  question of
orthogonality for more restrictive
class of functions $\svf s $; therein the tiles are rectangles in the phase plane. 
And the analysis shows
that orthogonality is linked to the boundedness of the the related maximal function.

There is a second form of this Theorem in which a multiscale object plays the role of a 
single tile.
We call that object a {\sl cluster}.  

\bthm d.cluster  
Let $\{\svf s \mid s\in\C\}$ be adapted to a set of tiles $\C$.  We call $\{\svf s \mid 
s\in\C\}$
a {\sl cluster with shadow } $I$  if the following four conditions are met. $(a)$ $\{s\in\C\}$ are pairwise
disjoint, $(b)$ for all $s\in\C$, $\sw s \subset I$,
\md4
(c)\ s\not=s'\in\C,\ \szw s \cap\szw s' \not=\emptyset\quad\text{implies}\quad 
\ip \svf s ,\svf s' .=0,
\\
(d)\ \sum_{s\in\C}\abs{\sw s }^{-1}\Biggl[\int_{I^c}(1+\zs(x,\sw s ))^{-d-1}\;dx\biggr]^2
{}\le{}C_1\abs{ I}.
\endmd
\ethm

\bthm t.cluster  Let $\S$ be a collection of tiles which is a disjoint union of 
subcollections 
$\{\C_I\mid I\in\I\}$. Let  $\{\svf s \mid s\in\S\}$ be adapted to $\S$, so that for each $I\in\I$, 
$\{ \svf s \mid s\in\C_I\}$ is a cluster. Suppose that \e.grid/ and \e.geo/ hold.  Finally, suppose
that the clusters $\C_I$ satisfy this condition.
\md1 
\Label e.cl.disjoint 
s\not=s'\in\S,\ s\in\C_I,\ \szw s \subset\szw s' 
\quad\text{implies}\quad s'\in\C_I\text{${}$ or ${}$}\sw s' \cap I=\emptyset.
\endmd
 Define 
\md0 
SQ(I,f)^2:=\frac1{\abs{I}} \sum_{s\in\C_I}\abs{\ip f,\svf s .}^2,\quad
I\in\I.
\endmd
Then, under the assumptions $SQ(I,f)\ge\zl$ for all $I\in\I$ and
 $\abs{\ip f,\svf s .}\le2\zl\sqw s $ for all $s\in\S$, we have 
 \md0
\sum_{I\in\I}\abs{I}\le{}\zl^{-2}K(1+ C_1^{1/2}C_0^2)\lVert
f\rVert_2^2.
\endmd
\ethm

The non--linear form of the hypotheses of this last Lemma preclude a natural
formulation of a weak--type inequality.

 These Theorems can also be used to study two complicated operators of harmonic analysis, 
namely
 Carleson's operator controlling the maximum partial Fourier sums \cite{carleson} and the 
bilinear
 Hilbert transform \cite{lt.pnas}.  A result clearly related to these Theorems can be found 
in a
 neglected paper of Prestini \cite{prestini}.  But the first forms of these Theorems, again 
for
 special tiles, is in \cite{lacey}.

\section{Proofs of Theorems}

 We will need a precise estimate of $\ip \svf s ,\svf s' .$, which is the purpose of 
 
 \bthm  l.ip    Let $s$ and $s'$ be two  tiles with $\sw s' -c(\sw s' )\subset\sw s -c(\sw s )$ but $\sw
 s' \cap\sw s  =\emptyset$.
  Let $\{\svf s ,\svf s' \}$ be 
 adapted to $\{s,s'\}$.  Then there is a constant $K$ so that 
 \md1\Label e.ss' 
 \abs{\ip \svf s ,\svf s' .}\le{}K{C_0^2}\sqrt{\frac{\abs{\sw s' }}{\abs{\sw s 
 }}}
 \inf_{x\in \sw s' }(1+\zs(x,\sw s ))^{-d-5}.
 \endmd
 \ethm
 
 \begin{proof} Heuristically, the Lemma follows from the estimate 
 \md0
 \abs{\ip \svf s ,\svf s' .}\simeq\int_{\sw s' }\abs{\svf s (x)\svf s' (x)}\;dx\le{}
 \norm \svf s .L^\infty(\sw s' ).\norm \svf s' .1.
 \endmd
 and \e.space/.  Of course the first step  must be made precise. \smallbreak
 
  Observe that 
\md0
\abs{\svf s (x)}\le{}K\frac{C_0}{\sqw s }\int_0^\zI \ind {a\sw s 
}(x)\;\frac{da}{(1+a)^{2d+6}},
\endmd
which is a consequence of \e.space/. We can estimate 
\md2
\abs{\ip \svf s ,\svf s' .}\le{}&
    K\frac{C_0^2}{\sqrt{\abs{\sw s }\abs{\sw s' }}}\int_0^\zI\int_0^\zI  \ip \ind {a\sw s 
},\ind {\za\sw s' }.\;
    \frac{da}{(1+a)^{2d+6}}\frac{d\za}{(1+\za)^{2d+6}}
    \\
    {}\le{}& KC_0^2\sqrt{\frac{\abs{\sw s' }}{\abs{\sw s 
 }}}\int_0^\zI\int_0^\zI\ind {\{a\sw s \cap\za\sw s' \not=\emptyset\}}
 \frac{da}{(1+a)^{2d+6}}\frac{\za^dd\za}{(1+\za)^{2d+6}}.
    \endmd
 The point to exploit is that $a\sw s \cap\za\sw s' =\emptyset$ for $a,\za\le{}\z1\zs:=\sup\{ a\mid
 aW_s\cap aW_{s'}\not=\emptyset\}$.

 We bound the double integral above by breaking the region of integration into four pieces.
  Take 
$R^0=[0,\z1\zs)$,
    $R^1=[\z1\zs,\zI)$ and define regions in the $(a,\za)$ plane by 
    \md0
    R_{\zve_1\zve_2}=R^{\zve_1}\times R^{\zve_2},\qquad \zve_i\in\{0,1\}.
    \endmd
   The integral over $R_{00}$ is zero by the choice of $\z1\zs$.  The integral over $R_{01}$ 
is  
    \md0
     \int_0^{\z1\zs}\int_{\z1\zs}^\zI 
\frac{da}{(1+a)^{2d+6}}\frac{\za^d d\za}{(1+\za)^{2d+6}}
     {}\le{} K (1+\z1\zs )^{-d-5} .
     \endmd
    We have a similar estimate for the integral over $R_{10}$.  Finally, the 
estimate over
     $R_{11}$ is    
\md0
     \int_{\z1\zs}^\zI\int_{\z1\zs}^\zI  
   \;
    \frac{da}{(1+a)^{2d+6}}\frac{\za^d d\za}{(1+\za)^{2d+6}}
\le{} K (1+\z1\zs)^{-d-5}.
         \endmd
    These estimates supply us with 
    \md0
    \abs{\ip \svf s ,\svf s' .}\le{} K C_0^2 \sqrt{\frac{\abs{ \sw s' }}{\abs{ \sw s }}} 
(1+\z1\zs)^{-d-5},
 \endmd
 which  is quite close to our claim.
    \smallskip

    To finish, observe that $\z1\zs\ge1$ as $\sw s $ and $\sw s' $ are disjoint.
     We shall also see that $\zs(x,\sw s )\le1+2\z1\zs$ for all $x\in\sw
    s' $.  These two points finish the proof of the Lemma.
    
      Indeed, we can assume that $\sw s $ and $\sw s' $ are closed sets.  Let $\z1x$ be a point in 
    $\z1\zs\sw s \cap\z1\zs\sw s' $.  Note
that for any point $x\in\sw s' $, 
\md0
x-c(\sw s )=(x-c(\sw s' ))+(c(\sw s' )-\z1x)+(\z1x -c(\sw s )).
\endmd
Of the three terms on the right, the first is in $\sw s' -c(\sw s' )\subset \sw s -c(\sw s )$, by
assumption,  the second is in $\z1\zs(\sw s' -c(\sw s' ))\subset\z1\zs(\sw s -c(\sw s ))$ by symmetry of
$\sw s' $ and assumption, 
and the third is in $\z1\zs(\sw s -c(\sw s ))$. The convexity of $\sw s $ then implies that 
$x\in(2\z1\zs+1)\sw s $.  Hence, $\zs(x,\sw s )\le1+2\z1\zs$ for all $x\in\sw
    s' $ as was to be shown. 
    \end{proof}
    
\begin{proof}[Proof of \t.count1/]
It is sufficient to prove a different assertion.  For   $f\in L^2(\ZR^d)$,
we assume that 
\md1 \Label e.b 
1\le{}\frac{\abs{\ip f,\svf s .}}{\sqw s }\le{}2,\quad s\in\S,
\endmd
and prove that 
\md1\Label e.-2 
\sum_{s\in\S}\abs{\sw s }\le{}(1+KC_0^2)\lVert f\rVert_2^{2}.
\endmd
For $k\ge0$ define $\S_k:=\{s\in\S\mid 2^{k}\le{}\abs{\sw s }^{-1/2}\abs{\ip  f,\svf s
.}\le{}2^{k+1}\}$.  One sees that the sum over this set of tiles is at most $2^{-2k}$ times 
the upper
bound in \e.-2/.  The Theorem is then  established in the case of $\zl=1$, but this
is sufficient as $f\in L^2(\ZR^d)$ is arbitrary.

For the proof of \e.-2/ we can assume that $\S$ is a finite collection, so that {\sl a 
priori} the
quantity
\md0
B:=\NORm \sum_{s\in\S}\langle f,\svf s \rangle \svf s .2.
\endmd
is finite.  It suffices to estimate $B$, for by using \e.b/
we see that 
\md1\Label e.W 
\sum_{s\in\S}\abs{\sw s }\le\sum_{s\in\S}\abs{\ip f,\svf s .}^2\; =\; {}
\IP f,\sum_{s\in\S}\langle f,\svf s \rangle \svf s .\le{}B\norm f.2..
\endmd

  We expand $B^2$ into diagonal and off--diagonal terms. 
\md1\Label e.B2  
B^2=\Bigl\lVert \sum_{s\in\S}\langle f,\svf s \rangle \svf s 
\Bigr\rVert_2^2\le{}\sum_{s\in\S}\abs{\langle f,\svf s \rangle  }^2+{}
2\sum_{s\in\S}|\langle f,\svf s \rangle |\z0O(s),
\endmd
where we define $\S(s)=\{s'\in\S-\{s\}\mid \szw s \subset\szw s' ,\langle \svf s
,\svf s' \rangle\not=0 \},$ and 
\md1\Label e.O  
\z0O(s):=\sum_{s'\in\S(s)}\abs{\langle \svf s ,\svf s' \rangle \ip \svf s' ,f.}.
\endmd
Recall that $\ip \svf s
,\svf s' .\not=0$ only if  $\szw s \cap \szw s' \not=\emptyset$.  But then from
the grid structure, we may assume that $\szw s \subset\szw s' $.

We have already seen that the diagonal term is dominated by $B\norm f.2.$, so that the term $\z0O(s)$ 
is our
concern.  Using \e.ss'/ and the upper bound on $\ip \svf s' ,f.$ we have 
\md2
\z0O(s)\le{}&KC_0^2\abs{\sw s }^{-1/2}\sum_{s'\in\S(s)} \inf_{x\in\sw s' 
}(1+\zs(x,\sw s
))^{-d-1}\abs{\sw s' }
\\
{}\le{}&KC_0^2\abs{\sw s }^{-1/2}\int_{(\sw s )^c}(1+\zs(x,\sw s
))^{-d-1}\;dx
\\
{}\le{}&KC_0^2\abs{\sw s }^{1/2}.
\endmd
For the middle line above,  the sets 
$\szw s' $
for $s'\in\S(s)$  contain $\szw s $.  But the tiles are disjoint, thus the sets $\sw s' $ are pairwise disjoint and contained in 
$(\sw s )^c$.

Therefore, the off diagonal term is, by \e.b/ and \e.W/,
\md0
\sum_{s\in\S}|\langle f,\svf s \rangle |\z0O(s)\le{}
  KC_0^2 \sum_{s\in\S} \abs{\sw s } 
  {}\le{}
  KC_0^2B\norm f.2..
  \endmd
  Combining this with \e.B2/ we see that 
  \md0
  B^2\le{}B\norm f.2.+ KC_0^2B\norm f.2.,
  \endmd
  which gives the desired upper bound $B$.
  \end{proof}

\begin{proof}[Proof of \t.cluster/.]
It suffices to consider the case of $\zl=1$.   The initial steps are just as before.  We assume that the collection of tiles 
is finite and  set 
\md0
B:=\NORm \sum_{s\in\S}\langle f,\svf s \rangle \svf s .2.
\endmd
and estimate $B$. This is sufficient since 
\md0
\sum_{I\in\I}\abs{I}\le{}\sum_{s\in\S}\abs{\ip f,\svf s .}^2\le{}B\norm f.2..
\endmd
 Then by Cauchy--Schwartz,
\md2
B^2\le{}&\sum_{s\in\S}\abs{\langle f,\svf s \rangle  }^2+{}
2\sum_{s\in\S}|\langle f,\svf s \rangle |\z0O(s),
\\
{}\le{}&B\norm f.2.+[B\norm f.2.]^{1/2}\biggl[\sum_{s\in\S}\abs{\z0O(s)}^2\biggr]^{1/2}
\endmd
Here as before, $\S(s)$
and $\z0O(s)$ are as in \e.O/.
 Then by \e.cl.disjoint/ and \d.cluster/ (c) the sets $\{\sw s' \mid s'\in\S(s)\}$ are 
 contained in $I^c$ if $s\in\C_I$.  But they are also pairwise disjoint.  Indeed, consider
 $s',s''\in\S(s)$.
 If they fall into the same cluster, they are disjoint by assumption on clusters.  Suppose they are 
  in different clusters.  Assuming as we may that 
 $\szw s' \subset \szw s'' $, we see that  \e.cl.disjoint/  implies 
 $\sw s' \cap \sw s'' =\emptyset$.
 
 Recalling \l.ip/ and  that we have  the upper bound $\abs{\ip f,\svf s'
 .}\le2\sqw s' $  by assumption, we see
 that  
 \md2
 \z0O(s)\le{}&
  KC_0^2\abs{\sw s }^{-1/2}\sum_{s'\in\S(s)}\inf_{x\in\sw s' 
}(1+\zs(x,\sw s ))^{-d-1}\abs{\sw s' }
\\
{}\le{}&
 KC_0^2\abs{\sw s }^{-1/2}\int_{I^c}(1+\zs(x,\sw s
))^{-d-1}\;dx
\endmd
if $s\in\C_I$.  Hence by the definition of a cluster
\md2
\sum_{s\in\C_I}\z0O(s)^2\le{}&
     KC_0^4\sum_{s\in\C_I}\abs{\sw s }^{-1}\Biggl[\int_{I^c}(1+\zs(x,\sw s 
))^{-d-1}\;
     {dx} \Biggr]^2
  \\
  {}\le{}&
     KC_1C_0^4\abs{I}.
     \endmd
 And so,  $B^2\le{}B\norm f.2.+KC_1^{1/2}C_0^2B\norm f.2.$, which proves the 
Theorem.
 \end{proof}    

\section{Counterexample}  We demonstrate the optimality of  the $L^{2,\zI}$ norm in our  Theorem.
It suffice to consider the first Theorem on $\ZR$.  And this we will do with the collection of disjoint
rectangles
\md0
\S:=\{[2^j,2^{j+1})\times[(n-1/2)2^{-j},(n+1/2)2^{-j})\mid j,n\in\ZZ\}.
\endmd
Let $\zvf$ denote a Schwartz function with $\zvf(x)>0$ for all $x$, $L^2$ norm one, and $\hat\zvf$ supported on
$[-1/2,1/2]$.  For $s=W_s\times\zW_s\in\S$ define 
\md0
\zvf_s(x):=e^{2\zp i c(\zW_s)x}\abs{W_s}^{-1/2}\zvf\Bigl(\frac{ x-c(W_s)}{\abs{W_s}}\Bigr).
\endmd
Clearly, $\{\zvf_s\mid s\in\S\} $ is adapted to $\S$.

Take $f=\ind{[-1,0)}$.  For integers $j\ge0$ and $\abs{n}<2^{j-1}$, 
let $s=[2^j,2^{j+1})\times[(n-1/2)2^{-j},(n+1/2)2^{-j})$.  We have 
\md2
2^{-j/2}\ip f,\zvf_s.=&2^{-j}\int_{-1}^0e^{2\zp i n2^{-j}x}\zvf\Bigl(\frac{ x-32^{j-1}}{2^j}\Bigr)\;dx
\\
{}={}&e^{3\zp in}\int_{-\frac{3}{2}-2^{-j}}^{-\frac{3}{2}}e^{2\zp i nx}\zvf(x)\;dx
\endmd
Since $\zvf(-3/2)>0$, we see that  $\abs{ 2^{-j/2}\ip f,\zvf_s.}\ge{}c2^{-j}$.  This estimate is uniform in
$n$ and $j$ as we have specified them.  Thus, for all $0<\zl<1$, 
\md0
\zl^2\sum_{s\in\S}\ind{}\{ \abs{\ip f,\zvf_s. } \ge\zl \sqrt{\abs{W_s}}\}\abs{W_s}\ge{}c'. 
\endmd
Thus, for any finite $t$,
\md0
\NOrm \Bigl\{\frac{\abs{\ip f,\svf s .}}{\sqrt{\abs{W_s}} }\ind {W_s} \mid s\in\S\Bigr\}
.L^{2,t}(\ZR^d\times \S).=\infty.
\endmd

\bigskip

\begin{tabular}{ll}
Jose Barrionuevo & \hskip .45in Michael Lacey \\
Department of Mathematics and Statistics & \hskip .45in School of Mathematics  \\
University of South Alabama &\hskip .45in Georgia Institute of Technology \\
Mobile AL 36688 &\hskip .45in Atlanta GA 30332  \\
   \\
\tt jose@mathstat.usouthal.edu &\hskip .45in\tt lacey@math.gatech.edu  \\
&\hskip .45in \tt http://www.math.gatech.edu/\~{}lacey \\
\end{tabular}


\begin{thebibliography}{99}
 \bibitem{carleson}L. Carleson. { ``On
convergence and growth of partial sums of
Fourier series."} {\it Acta Math.} {\bf 116} (1966) {pp. 135-157}.
 

\bibitem{lacey} M.T. Lacey.  { ``The bilinear Hilbert transform is pointwise finite."} {\it Rev.
Mat.} {\bf 13} (1997) 403---469.

\bibitem{lt} M.T. Lacey, C.M. Thiele.  { ``Bounds for the bilinear Hilbert transform on 
$L^p$." }  
Proc. Nat. Acad. Sci. {\bf 94} (1997) 33---35.

\bibitem{lt.pnas} M.T. Lacey, C.M. Thiele. {``$L^p$ Bounds for the bilinear Hilbert 
transform,  $p>2$."} {\sl Ann. Math.} {\bf 146} (1997) 693---724.

\bibitem{brush} F.G. Meyer, R.R. Coifman. ``Brushlets:  A tool for directional image analysis and
 image  compression."  {\sl Appl. Comp. Harmonic Anal.} {\bf 4} (1997) 188---221.

\bibitem{prestini}  E. Prestini.  { ``On the two proof of pointwise convergence of Fourier 
series."} 
{\it Amer. J. Math.} {\bf 104} (1982) 127---139.
 

\bibitem{w}  M. V. Wickerhauser. {\it Adapted Wavelet Analysis from Theory to Software} A K 
Peters Press,
1994.

\end{thebibliography}
\end{document}